\documentclass[a4paper]{article}
\usepackage{amsmath}
\usepackage{amssymb}
\usepackage{amsfonts}
\usepackage{theorem}
\usepackage{blkarray}
\usepackage{arydshln}
\usepackage{multirow}
\usepackage{graphicx}

\newcommand{\bigzerou}{%
\smash{\lower1.7ex\hbox{\bg 0}}}
\newtheorem{theorem}{Theorem}

\newtheorem{Rem}{Remark}
\newtheorem{lem}{Lemma}

\newcommand{\ba}{\begin{eqnarray}}
\newcommand{\ea}{\end{eqnarray}}
\newcommand{\no}{\nonumber}

\def\d{{\partial}}
\newcommand{\mapright}[1]{%
\smash{\mathop{%
\hbox to 1.0cm{\rightarrowfill}}\limits^{#1}}}
\newcommand{\mapleft}[1]{%
\smash{\mathop{%
\hbox to 1.3cm{\leftarrowfill}}\limits^{#1}}}
\allowdisplaybreaks
\begin{document}
\title{
\begin{flushright}
  \begin{minipage}[b]{5em}
    \normalsize
    ${}$      \\
  \end{minipage}
\end{flushright}
{\bf Hori's Equation for Gravitational Virtual Structure Constants of Calabi-Yau Hypersurface in $CP^{N-1}$}}
\author{Masao Jinzenji${}^{(1)}$, Kohki Matsuzaka${}^{(2)}$ \\
\\ 
${}^{(1)}$ \it Department of Mathematics,  \\
\it Okayama University \\
\it  Okayama, 700-8530, Japan\\
\\
${}^{(2)}$\it Institute for the Advancement of Graduate Education, \\
\it Hokkaido University \\
\it  Kita-ku, Sapporo, 060-0810, Japan\\
\\
\it e -mail address: 
\it\hspace{0.0cm}${}^{(1)}$ pcj70e4e@okayama-u.ac.jp \\
\it\hspace{3.5cm}${}^{(2)}$ kohki@math.sci.hokudai.ac.jp } 
\maketitle
\begin{abstract}
In this paper, we give a proof of Hori's equation for intersection numbers of the moduli space of quasimaps from $CP^{1}$ with (2+1) marked points to $CP^{N-1}$ by using localization technique.
\end{abstract}
\section{Introduction}

In this paper, we discuss the following two intersection numbers of the moduli spaces $\widetilde{Mp}_{0,2}(N,d)$ and $\widetilde{Mp}_{0,2|1}(N,d)$ of quasimaps from $CP^{1}$ to $CP^{N-1}$, which were constructed  
in the series of literatures of our group \cite{Jin1, Jin3, JS, S}.
The first one is, 
\begin{align}
&w \left( \sigma_{j} ({\cal O}_{h^{a}}) {\cal O}_{h^{b}} \right)_{0,d} \notag \\
&:= \int_{\widetilde{Mp}_{0,2}(N,d)} (\psi_{0})^{j} \cdot \mbox{ev}_{0}^{*} (h^{a}) \cdot \mbox{ev}_{\infty}^{*} (h^{b}) \cdot c_{top} ({\cal E}_{d}),
\end{align}
which was the main ingredient of our previous paper \cite{JM}. 

The second one is, 
\begin{align}
&w \left( \sigma_{j} ({\cal O}_{h^{a}}) {\cal O}_{h^{b}} \middle| {\cal O}_{h} \right)_{0,d} \notag \\
&:= \int_{\widetilde{Mp}_{0,2|1}(N,d)} (\psi^{\prime}_{0})^{j} \cdot \mbox{ev}_{0}^{*} (h^{a}) \cdot \mbox{ev}_{\infty}^{*} (h^{b}) \cdot \mbox{ev}_{1}^{*} (h) \cdot c_{top} ({\cal E}_{d}^{\prime}),
\end{align}
which we call ``gravitational virtual structure constant''. The reason why we call it gravitational virtual strucutre constant comes from the fact that our original virtual structure constant $\tilde{L}_{n}^{N,N,d}$ \cite{CJ} is expressed as $2+1$-pointed intersection number of $\widetilde{Mp}_{0,2|1}(N,d)$:
\ba
\tilde{L}_{n}^{N,N,d}=\frac{1}{N}w({\cal O}_{h^{N-2-n}}{\cal O}_{h^{n-1}}|{\cal O}_{h})_{0,d}.
\ea

In the above formulas, $\psi_{0}$ (resp. $\psi^{\prime}_{0}$) is Mumford-Morita-Miller class, i.e., the first Chern class of the line bundle on $\widetilde{Mp}_{0,2}(N,d)$ (resp. $\widetilde{Mp}_{0,2|1}(N,d)$) whose fiber at each point of the moduli space is given as cotangent space at the first marked point $0 \in CP^{1}$ or corresponding point in genus $0$ semi-stable curve used in compactification of  the moduli space.
${\cal E}_{d}$ (resp. ${\cal E}_{d}^{\prime}$) is orbi-bundle bundle on $\widetilde{Mp}_{0,2}(N,d)$ (resp. $\widetilde{Mp}_{0,2|1}(N,d)$) corresponding to the condition that images of quasimaps lie inside the Calabi-Yau hypersurface. Construction of these orbi-bundles will be discussed in Section 3. $\mbox{ev}_{0}$, $\mbox{ev}_{\infty}$ and $\mbox{ev}_{1}$ are the evalation maps at $0\in CP^{1}$, $\infty\in CP^{1}$ and $z_{1}\in CP^{1}-\{0,\infty\}$ respectively, and $h$ is the hyperplane class of $CP^{N-1}$.

We have to make several remarks on these intersection numbers. In \cite{JM}, the definition of  $w \left( \sigma_{j} ({\cal O}_{h^{a}}) {\cal O}_{h^{b}} \right)_{0,d}$ was given by using 
elements of Chow ring of $\widetilde{Mp}_{0,2}(N,d)$, whose explicit structure was given by Saito \cite{S}. On the other hand, in the above definition, we use expected  ``Mumford-Morita-Miller classes'' of the moduli space of quasimaps. Moreover, we haven't determined explicit structure of Chow ring of  $\widetilde{Mp}_{0,2|1}(N,d)$. Then ``how can we compute the above intersection numbers ?'' The key tool is the localization technique, which enables us to compute intersection numbers without knowing explicit structure of Chow ring.
As Givental did in \cite{giv}, the above intersection numbers can be computed explicitly by using information of moduli space of complex structure of weighted stable genus 0 curves
$\overline{M}_{0,2|1}$ \cite{AG,LM,MOP}. \footnote{In practice, $\overline{M}_{0,2|1}$ turns out to be just a point, and we don't need further information for computation of $w \left( \sigma_{j} ({\cal O}_{h^{a}}) {\cal O}_{h^{b}} \middle| {\cal O}_{h} \right)_{0,d}$. } We can show that the value of $w \left( \sigma_{j} ({\cal O}_{h^{a}}) {\cal O}_{h^{b}} \right)_{0,d}$ computed 
from the above definition and localization technique completely agrees with the one evaluated from the definition used in \cite{JM}. Of course, we can also compute  $w \left( \sigma_{j} ({\cal O}_{h^{a}}) {\cal O}_{h^{b}} \middle| {\cal O}_{h} \right)_{0,d}$ by combining Givental's discussion and our methods used in \cite{Jin1, Jin3, JS}.  

Using these results, we prove the following theorem.
\begin{theorem}{\bf (Main Theorem)}
For $j > 0$, the following equality holds.
\begin{equation}
w \left( \sigma_{j} ({\cal O}_{h^{a}}) {\cal O}_{h^{b}} \middle| {\cal O}_{h} \right)_{0,d} = d \cdot w \left( \sigma_{j} ({\cal O}_{h^{a}}) {\cal O}_{h^{b}} \right)_{0,d} + w \left( \sigma_{j-1} ({\cal O}_{h^{a+1}}) {\cal O}_{h^{b}} \right)_{0,d}.
\end{equation}
\label{main}
\end{theorem}
This theorem corresponds to quasimap analogue of ``Hori's equation'' \cite{hori} for gravitational Gromov-Witten invariants in the case of one K\"ahler insertion.

On the other hand, we have proved the following theorem in \cite{JM}:
\begin{theorem}
\ba
\frac{d}{N}w(\sigma_{j}({\cal O}_{h^{N-2-j}}){\cal O}_{h^{-1}})_{0,2}+
\frac{1}{N}w(\sigma_{j-1}({\cal O}_{h^{N-1-j}}){\cal O}_{h^{-1}})_{0,2}=\frac{1}{j!}\frac{\d^{j}}{\d \epsilon^{j}}\left(\frac{\prod_{r=1}^{Nd}(r+N\epsilon)}{\prod_{r=1}^{d}(r+\epsilon)^N}\middle) \right|_{\epsilon=0}.
\ea
\label{pmain}
\end{theorem}  
We can easily see that Theorem \ref{main} is correct even when $b=-1$ (see Remark \ref{final} at the end of this paper). Hence we obtain,
\ba
\frac{1}{N}w \left( \sigma_{j} ({\cal O}_{h^{N-2-j}}) {\cal O}_{h^{-1}} \middle| {\cal O}_{h} \right)_{0,d}= \frac{1}{j!}\frac{\d^{j}}{\d \epsilon^{j}}\left(\frac{\prod_{r=1}^{Nd}(r+N\epsilon)}{\prod_{r=1}^{d}(r+\epsilon)^N}\middle) \right|_{\epsilon=0}.
\ea
In this way, we can provide concrete footing with the assertion that period integrals used in mirror computation of genus $0$ Gromov-Witten invariants of Calabi-Yau hypersurface in $CP^{N-1}$ and correponding Givental's $I$-function are given as generating functions of  gravitational virtual structure constants \cite{JM}.

This paper is organized as follows. In Section 2, we present outline of construction of the moduli space $\widetilde{Mp}_{0,2|1}(N,d)$ following the construction given in \cite{JS}.
In Section 3, we briefly discuss construction of the orbi-bundle ${\cal E}_{d}^{\prime}$ on  $\widetilde{Mp}_{0,2|1}(N,d)$, which was omitted in the discussion of \cite{JS}.  
In Section 4, we prove the main thorem by using localization technique and results given in \cite{Jin1,Jin3,JS}.

\vspace{2em}

{\bf Acknowledgment} 
We would like to thank Prof.~G. Ishikawa and Prof.~A.Tsuchida for kind encouragement.  Our research is partially supported by JSPS grant No. 17K05214.  

\section{The Moduli Space $\widetilde{Mp}_{0,2|1}(N,d)$}

Following \cite{JS}, we review construction of the moduli space $\widetilde{Mp}_{0,2|1}(N,d)$ of quasimaps of degree $d$ from $CP^{1}$ with $2+1$ marked points to $CP^{N-1}$. 
First, we construct bulk part $Mp_{0,2|1}(N,d)$ of the moduli space.  
A degree $d$ \textit{quasimap} $\varphi (s:t)$ from $CP^{1}$ in  $Mp_{0,2|1}(N,d)$ is given by using a vector valued polynomial
$ \varphi(s,t)=\sum_{i =0}^{d} \mathbf{a}_{i} s^{d-i} t^{i}\;\;(\mathbf{a}_{i}\in \mathbb{C}^{N})$ as follows: 
\begin{equation}
\varphi (s:t) :=[\varphi(s,t)], 
\end{equation}
where $[*]$ means taking projective equivalence class of $*$. We impose the condition $\mathbf{a}_{0}, \mathbf{a}_{d}\neq {\bf 0}$ so that the
images $\varphi(1:0)$ and $\varphi(0:1)$ are well-defined in $CP^{N-1}$.
From the above definition of $\varphi(s:t)$, we can easily see  that it allows points in $CP^{1}$ other than $0, \infty$
whose images cannot be defined. Explicitly, we can factorize $\varphi(s,t)$ into the form, 
\begin{equation}
\varphi(s,t)=\sum_{i =0}^{d} \mathbf{a}_{i} s^{d-i} t^{i}= p(s,t) \cdot\biggl( \sum_{k = 0}^{d- \mathrm{deg}(p)} \mathbf{b}_{k} s^{d- \mathrm{deg}(p) -k} t^{k}\biggr),
\label{factor}
\end{equation}
so that $[ \sum_{k = 0}^{d- \mathrm{deg}(p)} \mathbf{b}_{k} s^{d- \mathrm{deg}(p) -k} t^{k}]$ defines a well-defind degree $d-\mbox{deg}(p)$ map from $CP^{1}$ to $CP^{N-1}$.
Then image of $\varphi$ cannot be defined at $(\alpha:\beta)\in CP^{N-1}$ that satisfies $p(\alpha,\beta)=0$. 
  
Then we introdude $2+1$ marked points $(0, \infty | z_1)$ in $CP^{1}$. Here $0=(1:0), \infty=(0:1)$ and $z_{1}$ is a point $(1:z_{1})\in CP^{1}$ which do not coincide 
with $0$ and $\infty$. Therefore, we can label $2+1$ marked points by $z_{1}\in CP^{1}-\{0,\infty\}=\mathbb{C}^{\times}$. Then we can represent a pair of $\varphi(s,t)=\sum_{i =0}^{d} \mathbf{a}_{i} s^{d-i} t^{i}$ and a set of $2+1$
marked points  $(0, \infty | z_1)$ as follows:
\ba
(\mathbf{a}_{0},\mathbf{a}_{1},\cdots,\mathbf{a}_{d},z_{1})\in (\mathbb{C}^{N}-\{\mathbf{0}\})\times\mathbb{C}^{N}\times\cdots\times\mathbb{C}^{N}\times(\mathbb{C}^{N}-\{\mathbf{0}\})\times \mathbb{C}^{\times}.
\ea
Next, we introduce $\mathbb{C}^{\times}\times \mathbb{C}^{\times}$ action on the above set,
\ba
(\mathbf{a}_{0} , \dots , \mathbf{a}_{i}, \dots , \mathbf{a}_{d}, z_1 ) \rightarrow  
(\mu \mathbf{a}_{0} , \dots , \mu \nu^{i} \mathbf{a}_{i}, \dots , \mu \nu ^{d} \mathbf{a}_{d}, \nu^{-1} z_1)\;\; (\mu,\nu\in \mathbb{C}^{\times}).
\ea
With these set-up's, $Mp_{0,2|1}(N,d)$ is defined as the following quotient space:
\begin{equation}
Mp_{0,2|1} (N,d) := \biggl((\mathbb{C}^{N}-\{\mathbf{0}\})\times\mathbb{C}^{N}\times\cdots\times\mathbb{C}^{N}\times(\mathbb{C}^{N}-\{\mathbf{0}\})\times \mathbb{C}^{\times}\biggr) / ( \mathbb{C}^{\times}\times \mathbb{C}^{\times}),
\end{equation}
Let us denote equivalence class of $(\mathbf{a}_{0} , \dots , \mathbf{a}_{i}, \dots , \mathbf{a}_{d}, z_1 )$ under the $\mathbb{C}^{\times}\times \mathbb{C}^{\times}$ action by  
$[(\mathbf{a}_{0} , \dots , \mathbf{a}_{i}, \dots , \mathbf{a}_{d}, z_1 )]$. At this stage, we define evaluation maps $\mbox{ev}_{0}: Mp_{0,2|1} (N,d)\rightarrow CP^{N-1}$
and $\mbox{ev}_{\infty}: Mp_{0,2|1} (N,d)\rightarrow CP^{N-1}$ as follows.
\ba
&&\mbox{ev}_{0}([(\mathbf{a}_{0} , \dots , \mathbf{a}_{i}, \dots , \mathbf{a}_{d}, z_1 )]):=[\varphi(1:0)]=[\mathbf{a}_{0}],\no\\
&&\mbox{ev}_{\infty}([(\mathbf{a}_{0} , \dots , \mathbf{a}_{i}, \dots , \mathbf{a}_{d}, z_1 )]):=[\varphi(0:1)]=[\mathbf{a}_{d}].
\ea
We also define $\mbox{ev}_{1}: Mp_{0,2|1} (N,d)\rightarrow CP^{N-1}$, i.e., evaluation map at $z_{1}$, but in this case, we have to consider subtle point that image of the quasimap
at $z_{1}$ may not be defined. In order to avoid this subtlety, we use the unique factorization (\ref{factor}) and define $\mbox{ev}_{1}$ as follows.
\ba
\mbox{ev}_{1}([(\mathbf{a}_{0} , \dots , \mathbf{a}_{i}, \dots , \mathbf{a}_{d}, z_1 )]):=[\sum_{k = 0}^{d- \mathrm{deg}(p)} \mathbf{b}_{k}(z_{1})^{k}],
\ea  
We leave confirmation of well-definedness of the above definitions to the readers as an exercise. 

In order to compactify $Mp_{0,2|1} (N,d)$, we introduce genus $0$ semi-stable curves whose topological type are given by line graphs (or chain type) and consider chain of quasimaps.
Each quasimap has a $CP^{1}$ componet as its domein and two quasimaps from two $CP^{1}$'s connected by a nodal singularity are glued at the singularity.  

First, let us explain the process of compactification in the case of $\widetilde{Mp}_{0,2}(N,d)$, which is the moduli space of degree $d$ quasimaps with no marked points except for $0$ and $\infty$ \cite{Jin1,Jin3}. The bulk part $Mp_{0,2}(N,d)$ is defined as follows.
\begin{equation}
Mp_{0,2} (N,d) := \biggl((\mathbb{C}^{N}-\{\mathbf{0}\})\times\mathbb{C}^{N}\times\cdots\times\mathbb{C}^{N}\times(\mathbb{C}^{N}-\{\mathbf{0}\})\biggr) / ( \mathbb{C}^{\times}\times \mathbb{C}^{\times}),
\end{equation}
The $\mathbb{C}^{\times}\times \mathbb{C}^{\times}$ is given by, 
\ba
(\mathbf{a}_{0} , \dots , \mathbf{a}_{i}, \dots , \mathbf{a}_{d})\rightarrow  
(\mu \mathbf{a}_{0} , \dots , \mu \nu^{i} \mathbf{a}_{i}, \dots , \mu \nu ^{d} \mathbf{a}_{d})\;\; (\mu,\nu\in \mathbb{C}^{\times}).
\ea
Construction of $\widetilde{Mp}_{0,2}(N,d)$, or in other words, compactification of $Mp_{0,2} (N,d)$ is done by introducing stable curve
\begin{equation}
\left( \coprod_{i = 1}^{l} B_i \middle) \right/ \sim
\end{equation}
for some $l \ge 1$. In the above formula, each $B_i$ is a copy of $CP^1$ and $\sim$ means equivalence relations given by,
\begin{equation}
\infty_1 \sim 0_2, \ \infty_2 \sim 0_3, \ \dots , \ \infty_{l-1} \sim 0_l ,
\end{equation}
where $0_i (= (1:0))$ and $\infty_i (= (0:1))$ are two marked points of $B_i = CP^1$. 
Then we introduce chain of quasimaps $(\varphi_{1},\varphi_{2},\cdots,\varphi_{l})$ that satisfy the following conditions.
\begin{itemize}
\item[(i)] { $d_{i}$, the degree of $\varphi_{i}$, is a positive integer.}
\item[(ii)] { $\varphi_{i}(\infty_{i})=\varphi_{i+1}(0_{i+1})\;\;(i=1,\cdots,l-1).$} 
\end{itemize}
Then $\widetilde{Mp}_{0,2}(N,d)$ is constructed as follows. 
\begin{equation}
\widetilde{Mp}_{0,2}(N,d) := \coprod_{(d_1 , \dots , d_l ) \in S(d)} Mp_{0,2} (N,d_1 ) \times _{CP^{N-1}} \dots \times _{CP^{N-1}} Mp_{0,2} (N,d_l ) ,
\end{equation}
where $S(d)$ is set of orderd partitions of $d$:
\begin{equation}
S(d) := \{ (d_1 , \dots , d_l ) \; | \; 1\leq l\leq d,  d_1 , \dots , d_l \in \mathbb{Z}_{>0} , d_1 + \dots + d_l = d \;\} ,
\end{equation}
and $Mp_{0,2} (N,d_i ) \times _{CP^{N-1}} Mp_{0,2} (N,d_{i+1} )$ is fiber product via two evaluation maps:
\begin{equation}
\mbox{ev}_{\infty_i} : Mp_{0,2} (N,d_i ) \to CP^{N-1} \; \; (\varphi_i \mapsto \varphi_i (\infty_i )),
\end{equation}
and
\begin{equation}
\mbox{ev}_{0_{i+1}} : Mp_{0,2} (N,d_{i+1} ) \to CP^{N-1} \;  \; (\varphi_{i+1} \mapsto \varphi_{i+1} (0_{i+1} )).
\end{equation}
The fiber products in the above corresponds to the condition (ii) of the chain of quasimaps. 

As was mentioned in \cite{Jin1, S}, $\widetilde{Mp}_{0,2}(N,d)$ is also constructed as a quotient space (toric variety):
\begin{align}
&\widetilde{Mp}_{0,2}(N,d) \notag \\
&= \{ (\mathbf{a}_{0} , \dots , \mathbf{a}_{d} , u_1 , \dots , u_{d-1} ) \in \mathbb{C}^{N(d+1) + (d-1)} \; | \; \mathbf{a}_{0} , (\mathbf{a}_1 , u_1 ) , \dots , (\mathbf{a}_{d-1} , d_{d-1}) , \mathbf{a}_{d} \neq \mathbf{0} \} / (\mathbb{C}^{\times})^{d+1} , \label{qrep}
\end{align}
where $(d+1)$ $\mathbb{C}^{\times}$ actions are given by,
\begin{align} 
&(\mu_{0} , \dots , \mu_{d}) \cdot (\mathbf{a}_{0} , \dots , \mathbf{a}_{i} , \dots , \mathbf{a}_{d} , u_1 , \dots , u_{i-1} , u_i , u_{i+1} , \dots , u_{d-1}) \notag \\
&= (\mu_0 \mathbf{a}_{0} , \dots , \mu_i \mathbf{a}_{i} , \dots , \mu_{d} \mathbf{a}_{d} , \mu_0^{-1} \mu_1^2 \mu_2^{-1} u_1 , \dots , \mu_i^{-1} u_{i-1} , \mu_i^2 u_i , \mu_i^{-1} u_{i+1} , \dots , \mu_{d-2}^{-1} \mu_{d-1}^2 \mu_{d}^{-1} u_{d-1}) \notag\\
& \hspace{11cm} (\mu_0 , \dots , \mu_{d} \in \mathbb{C}^{\times}).
\end{align}
 New coordinates $u_1 , \dots , u_{d-1}$ represent the ``boundary locus'' of $\widetilde{Mp}_{0,2}(N,d)$. More precisely, the condition: 
\ba
&&u_{i_1} = \dots = u_{i_l} = 0 \quad (1 \le i_1 < \dots < i_l \le d-1 ), \notag \\
&&\hspace{2cm}u_{i} \neq 0 \quad (\mbox{otherwise}),
\ea
corresponds to the boundary component:
\begin{align}
Mp_{0,2} (N, i_1 ) &\times_{CP^{N-1}} Mp_{0,2} (N, i_2 - i_1 ) \times _{CP^{N-1}} \dots \notag \\
&\dots \times _{CP^{N-1}} Mp_{0,2} (N, i_{l} - i_{l-1} ) \times _{CP^{N-1}} Mp_{0,2} (N, d - i_{l} ).
\end{align}
In particular, a point $[ (\mathbf{a}_{0} , \dots , \mathbf{a}_{d} , u_1 , \dots , u_{d-1}) ]$ of $ \widetilde{Mp}_{0,2}(N,d)$ corresponds to a point of $Mp_{0,2} (N,d)$ if and only if $u_1 , \dots , u_{d-1} \neq 0$. 

In order to compactify $Mp_{0,2|1} (N,d)$, 
we only have to add information of position of $z_{1}$ to the above construction. 
Fundamentaly, $z_{1}$ belongs to $B_{i}-\{0_{i},\infty_{i}\}$ for some $i\in \{1,2,\cdots,l\}$,
but we have one subtlety that we have to describe the limit $z_1 \to 0_{i}$ or $z_{1} \to \infty_{i}$ when $z_{1}\in (B_{i}-\{0_{i},\infty_{i}\})$. 

 For example, as the case of stable maps, the situation $z_1 \to 0$ for the 2+1 marked points $\left(0,\infty \middle| z_1 \right)$ in $CP^{1}$ is expressed by adding ``infinitessimal component'' $A = CP^1$ with $2+1$ marked points $\left(0^{\prime} , \infty^{\prime}|z_1^{\prime}\right)$ in $A$. Then we consider, $$A\coprod CP^{1}/\sim,$$ (original $CP^{1}$ has now only two marked points $(0,\infty)$) where equivalence relation $\sim$ is given by $0 \sim \infty^{\prime}$. Then we can describe the situation $z_1 \to 0$ by using this semi-stable curve. 
The case $z_1 \to \infty$ is discribed in a simillar way. In order to proceed to the level of quasimaps, we impose the condition that the infinitessimal component $A$ is mapped to a point 
in $CP^{N-1}$, or  in other words, quasimap assigned to this component has degree $0$. Then the degree of freedom of the moduli space coming from this component $A$ is given by
$\overline{M}_{0,2|1}$, ``the moduli space of complex structure of $CP^{1}$ with $2+1$ marked points'' \cite{AG,LM,MOP}, but this moduli space is known to be given as ``a point''
because we can always set $z_{1}=1$ by using $\mathbb{C}^{\times}$ acion (automorphism group of $CP^{1}$ that fix $0$ and $\infty$).   

With these considerations, $\widetilde{Mp}_{0,2|1}(N,d)$ is constructed as follows.
\begin{equation}
\widetilde{Mp}_{0,2|1}(N,d) := \coprod_{(d_1 , \dots , d_l ) \in S(d)} \left( \coprod_{i=1}^{l} \mathrm{Nor}(d_1 , \dots , d_l ; i) \right) \coprod \left( \coprod_{i=0}^{l} \mathrm{Lim}(d_1 , \dots , d_l ; i) \right) ,
\end{equation}
where $\mathrm{Nor}(d_1 , \dots , d_l ; i)$ is normal part that describes the case $z_{1}\in (B_{i}-\{0_{i},\infty_{i}\})$:
\begin{align}
\mathrm{Nor}(d_1 , \dots , d_l ; i) &:= Mp_{0,2} (N,d_1 ) \times_{CP^{N-1}} \dotsb \notag \\
&\quad \dots \times_{CP^{N-1}} Mp_{0,2} (N,d_{i-1} ) \times_{CP^{N-1}} Mp_{0,2|1} (N,d_i ) \times_{CP^{N-1}} Mp_{0,2} (N,d_{i+1} ) \times_{CP^{N-1}} \dotsb \notag \\
&\quad \dots \times_{CP^{N-1}} Mp_{0,2} (N,d_l ) \;\;(i=1,2,\cdots,l), 
\end{align}
and $\mathrm{Lim}(d_1 , \dots , d_l ; i)$ is limit part that describes the case $z_1 \to \infty_{i}$ or $z_{1} \to 0_{i+1}$:
\begin{align}
\mathrm{Lim}(d_1 , \dots , d_l ; i) &:= Mp_{0,2} (N,d_1 ) \times_{CP^{N-1}} \dotsb \notag \\
&\quad \dots \times_{CP^{N-1}} Mp_{0,2} (N,d_i ) \times_{CP^{N-1}} (CP^{N-1} \times \overline{M}_{0,2|1} ) \times_{CP^{N-1}} Mp_{0,2} (N,d_{i+1} ) \times_{CP^{N-1}} \dotsb \notag \\
&\quad \dots \times_{CP^{N-1}} Mp_{0,2} (N,d_l )\;\;(i=0,1,2,\cdots,l), 
\end{align}
where each fiber product is defined in the same way as the case of $Mp_{0,2} (N,d)$. In the limit part, $CP^{N-1} \times \overline{M}_{0,2|1}$ between $Mp_{0,2} (N,d_i)$ and 
 $Mp_{0,2} (N,d_{i+1})$ describes the situation that $z_{1}\in B_{i} \to \infty_{i}$ and $z_{1}\in B_{i+1} \to 0_{i+1}$ (note that $\infty_{i}=0_{i+1}$) and the $CP^{N-1}$ in the left entry gives the
 information of the image of $A$ by degree $0$ quasimap (constant map) in $CP^{N-1}$. Since $\overline{M}_{0,2|1}$ is just a point, we have the following set-theoretical equality:
 \begin{align}
 \mathrm{Lim}(d_1 , \dots , d_l ; i) =& Mp_{0,2} (N,d_1 ) \times_{CP^{N-1}} \dotsb \times_{CP^{N-1}} Mp_{0,2} (N,d_i ) \times_{CP^{N-1}} Mp_{0,2} (N,d_{i+1} ) \times_{CP^{N-1}} \dotsb \notag 
 \\
&\quad \dots \times_{CP^{N-1}} Mp_{0,2} (N,d_l ).
 \end{align}

In closing discussion of this section, we introduce local coordinate system on $Mp_{0,2|1} (N,d)$. Let $U_{i|k} \ (i, k = 1 , \dots , N)$ be an open subset of $Mp_{0,2|1} (N,d)$:
\begin{equation}
U_{i|k} := \left\{ (\mathbf{a}_0 , \mathbf{a}_1 , \dots , \mathbf{a}_{d-1} , \mathbf{a}_d , z_1 ) \in Mp_{0,2|1} (N,d) \; \middle| \; a_0^i , a_d^k \neq 0 \right\} ,
\end{equation}
where $\mathbf{a}_{\bullet} = (a_{\bullet}^{1} , \dots , a_{\bullet}^{N})$, and $\phi_{i|k} : U_{i|k} \to \{ \mathbb{C}^{N-1} \times (\mathbb{C}^{N})^d \times \mathbb{C}^{N-1} \times \mathbb{C}^{\times} \} / \mathbb{Z}_d$ be a map:
\begin{align}
\phi_{i|k} ([\mathbf{a}_0 , \mathbf{a}_1 , \dots , \mathbf{a}_{d-1} , \mathbf{a}_d , z_1]) 
&:= \Bigg( \left( \frac{a_{0}^{1}}{a_{0}^{i}} , \dots , \frac{a_{0}^{i-1}}{a_{0}^{i}} , \frac{a_{0}^{i+1}}{a_{0}^{i}} , \dots , \frac{a_{0}^{N}}{a_{0}^{i}} \right) \notag \\
&\hspace{25pt} , \left[ \frac{\mathbf{a}_1}{(a_0^i )^{\frac{d-1}{d}} (a_d^k )^{\frac{1}{d}} } , \dots , \frac{\mathbf{a}_m}{(a_0^i )^{\frac{d-m}{d}} (a_d^k )^{\frac{m}{d}} } , \dots , \frac{\mathbf{a}_{d-1}}{(a_0^i )^{\frac{1}{d}} (a_d^k )^{\frac{d-1}{d}} } \right] \notag \\
&\hspace{25pt} , \left( \frac{a_{d}^{1}}{a_{d}^{k}} , \dots , \frac{a_{d}^{k-1}}{a_{d}^{k}} , \frac{a_{d}^{k+1}}{a_{d}^{k}} , \dots , \frac{a_{d}^{N}}{a_{d}^{k}} \right) \notag \\
&\hspace{25pt} , \left[ \frac{(a_d^k )^{\frac{1}{d}}}{(a_0^i )^{\frac{1}{d}}} z_1 \right] \Bigg) ,
\end{align}
where finite cyclic group $\mathbb{Z}_d$ acts simultaneously on the 2nd and the 4th factor as follows:
\ba
&&\mu_d \cdot (\mathbf{a}_1 , \dots , \mathbf{a}_m , \dots , \mathbf{a}_{d-1}) = ( \mu_d \mathbf{a}_1 , \dots , (\mu_d )^m \mathbf{a}_m , \dots , (\mu_d )^{d-1} \mathbf{a}_{d-1}) , \quad (\text{2nd}) \notag \\
&&\hspace{3.05cm}\mu_d \cdot z_1 = \mu_d z_1,  \quad (\text{4th}) \notag\\
&&\hspace{3.25cm}\left( \mu_d := \mathrm{exp} \left( \frac{2 \pi \sqrt{-1}}{d} \right) \right) . 
\ea
It is easy to see that $\phi_{i|k}$ is bijective, and $\{ ( U_{i|k} , \phi_{i|k} ) \}_{i,k}$ is a local chart on $Mp_{0,2|1} (N,d)$ (but the second factor has an orbifolod singularity
$\mathbf{a}_{1}=\mathbf{a}_{2}=\cdots=\mathbf{a}_{d-1}=\mathbf{0}$). Therefore, $Mp_{0,2|1} (N,d)$ is an orbifold. This situation is the same as the one of  $Mp_{0,2} (N,d)$ \cite{Jin1,Jin3}.

\section{The Orbi-Bundle ${\cal E}_{d}^{\prime}$}

In this section, we construct the orbi-bundle ${\cal E}_{d}^{\prime}$ on the moduli space $Mp_{0,2|1}(N,d)$ and then we extend it to the whole parts of $\widetilde{Mp}_{0,2|1}(N,d)$. In the same way as the case of Gromov-Witten invariants, vanishing locus of a global section of this orbi-bundle expresses subset of quasimaps whose images lie inside Calabi-Yau hypersurface $M$ in $CP^{N-1}$. Assume that $M$ is defined by degree $N$ homogeneous polynomial $F(X_1 , \dots , X_N )$ (e.g. $(X_1 )^{N} + \dots + (X_N )^{N}$) and let $f_{i} (\mathbf{a}_{0} , \dots , \mathbf{a}_{N})$ ($i = 0 , \dots , Nd$) be a degree $N$ homogeneous polynomial in $a_{\bullet}^{i}$ obtained from the following equality:
\begin{equation}
F \left( \sum_{i=0}^{d} a_{i}^{1} s^{d-i} t^{i} , \dots , \sum_{i=0}^{d} a_{i}^{N} s^{d-i} t^{i} \right) = \sum_{i=0}^{Nd} f_{i} (\mathbf{a}_{1} , \dots , \mathbf{a}_{N}) s^{Nd-i} t^{i} .
\end{equation}
The two $\mathbb{C}^{\times}$ actions used in defining $Mp_{0,2|1}(N,d)$ act on these homogeneous polynomials as follows.:
\begin{equation}
f_m (\mu \mathbf{a}_{0} , \dots , \mu \nu^i \mathbf{a}_i , \dots , \mu \nu^{d} \mathbf{a}_{N}) = \mu^N \nu^{m} f_m (\mathbf{a}_{0} , \dots , \mathbf{a}_{N}) \quad (m = 0 , \dots , Nd). \label{glosec}
\end{equation}
The condition ``$\mathrm{Im} \; \varphi \subset M$'' is equivalent to
\begin{equation}
f_{0} (\mathbf{a}_{0} , \dots , \mathbf{a}_{N}) = \dots = f_{Nd} (\mathbf{a}_{0} , \dots , \mathbf{a}_{N}) = 0,
\end{equation}
and (\ref{glosec}) tells us that this condition defines a well-defined subset of $Mp_{0,2|1}(N,d)$.
Then we construct the orbi-bundle ${\cal E}_{d}^{\prime}$ so that it has 
\begin{equation}
( f_{0} (\mathbf{a}_{0} , \dots , \mathbf{a}_{N}) , \dots , f_{Nd} (\mathbf{a}_{0} , \dots , \mathbf{a}_{N})).
\end{equation}
as a global section. 
Let $\{ ( U_{i|k} , \phi_{i|k} ) \}_{i,k}$ be the local chart on $Mp_{0,2|1} (N,d)$ introduced in the previous section. We define system of transition functions $\{ g_{ij|kl}^{(m)} \}_{(i,k),(j,l)}$ ($m = 0 , \dots , Nd$) on $\{ U_{i|k} \}_{i,k}$ as
\begin{equation}
g_{ij|kl}^{(m)} ([ \mathbf{a}_{0} , \dots , \mathbf{a}_{N} , z_1 ]) := \left( \frac{a_0^j}{a_0^i} \right)^{\frac{Nd-m}{d}} \left( \frac{a_d^l}{a_d^k} \right)^{\frac{m}{d}} \quad (\text{on} \ \ U_{ij|kl} := U_{i|k} \cap U_{j|l}).
\end{equation}
Then it gives a rank $1$ orbi-bundle ${\cal E}_{d}^{\prime(m)}$ which is isomorphic to 
\begin{equation}
{\cal O}_{(CP^{N-1})_0} \left( \frac{Nd-m}{d} \right) \otimes {\cal O}_{(CP^{N-1})_d} \left( \frac{m}{d} \right) := \mbox{ev}_{0}^{*} {\cal O}_{CP^{N-1}} \left( \frac{Nd-m}{d} \right) \otimes \mbox{ev}_{\infty}^{*} {\cal O}_{CP^{N-1}} \left( \frac{m}{d} \right).
\end{equation}
In the above, we denote by $(CP^{N-1})_{0}$ and $(CP^{N-1})_{d}$ the images of $\mbox{ev}_0$ and $\mbox{ev}_{\infty}$, respectively. By the definition of coordinate system $\{ ( U_{i|k} , \phi_{i|k} ) \}_{i,k}$ and the relation (\ref{glosec}), it is easy to see that a system of functions $\{ (s_{i|k}^{(m)} , U_{i|k}) \}_{i,k}$ defined by
\begin{equation}
s_{i|k}^{(m)} ([ \mathbf{a}_{0} , \dots , \mathbf{a}_{N} , z_1 ]) := \frac{f_{m} (\mathbf{a}_{0} , \dots , \mathbf{a}_{N})}{(a_0^i )^{\frac{Nd-m}{d}} (a_{d}^k )^{\frac{m}{d}}}
\end{equation}
gives a global section of ${\cal E}_{d}^{\prime(m)}$ that correponds to $f_m (\mathbf{a}_{0} , \dots , \mathbf{a}_{N})$. Hence we define ${\cal E}^{\prime}_{d}$ as follows.
\begin{align}
{\cal E}_{d}^{\prime} &:= \bigoplus_{m=0}^{Nd} {\cal E}_{d}^{\prime(m)} \notag \\
&\; \simeq \bigoplus_{m=0}^{Nd} \left( {\cal O}_{(CP^{N-1})_0} \left( \frac{Nd-m}{d} \right) \otimes {\cal O}_{(CP^{N-1})_d} \left( \frac{m}{d} \right) \right) . 
\end{align}
This discussion is the same as construction of the orbi-bundle ${\cal E}_{d}$ on $Mp_{0,2}(N,d)$ given in \cite{Jin1, Jin3} because $z_{1}$ plays no role on the condition that image of a quasimap lies inside the hypersurface.

Then we extend the above ${\cal E}_{d}^{\prime}$ to whole parts of $\widetilde{Mp}_{0,2|1}(N,d)$. We discuss separately the following two parts; (i) normal part $\mathrm{Nor}(d_ 1, \dots , d_l ; i)$ and (ii) limit part $\mathrm{Lim}(d_1 , \dots , d_l ; i)$.

In the case of the normal part, we first consider simpler case $Mp_{0,2|1} (N,d_1 ) \times_{CP^{N-1}} Mp_{0,2|0} (N,d_2 )\;\;(d=d_{1}+d_{2})$. A point in $Mp_{0,2|1} (N,d_1 ) \times_{CP^{N-1}} Mp_{0,2|0} (N,d_2 )$ is given as a chain of quasimaps $(\varphi_{1},\varphi_{2})$ that represents a stable quasimap defined on semi-stable curve obtained from gluing two copy of $CP^{1}$, 
$B_{1}$ and $B_{2}$,  by the equivalence relation $\infty_1 \sim 0_2$.
Then $\varphi_{1}$ and $\varphi_{2}$ must satisfy $\varphi_1 (\infty_1 ) = \varphi_2 (0_2 )$. Therefore, in order to extend ${\cal E}_{d}^{\prime}$ to $Mp_{0,2|1} (N,d_1 ) \times_{CP^{N-1}} Mp_{0,2|0} (N,d_2 )$, we have to glue ${\cal E}_{d_1}^{\prime}$ on $Mp_{0,2|1} (N,d_1 )$ and ${\cal E}_{d_2}$ on $Mp_{0,2}(N,d_{2})$ by the following exact sequence:
\begin{equation}
0 \to {\cal E}_{d}^{\prime}|_{Mp_{0,2|1} (N,d_1 ) \times_{CP^{N-1}} Mp_{0,2} (N,d_2 )} \to {\cal E}_{d_1}^{\prime} \oplus {\cal E}_{d_2} \to p_{1,2} ^{*} {\cal O}_{CP^{N-1}} (N) \to 0,
\end{equation}
where $p_{1,2}$ is a projection map:
\begin{equation}
p_{1,2} : Mp_{0,2|1} (N,d_1 ) \times_{CP^{N-1}} Mp_{0,2|0} (N,d_2 ) \to CP^{N-1} ; \; (\varphi_1 , \varphi_2 ) \mapsto \varphi_1 (\infty _1 ) \; (= \varphi_2 (0_2 )) .
\end{equation}
In general, in order to extend the orbi-bundle ${\cal E}_{d}^{\prime}$ to $\mathrm{Nor} (d_1 , \dots , d_l ; i)$, we use the  following exact sequence:
\begin{equation}
0 \to {\cal E}_{d}^{\prime}|_{\mathrm{Nor} (d_1 , \dots , d_l ; i)} \to \left( \bigoplus_{k=1}^{i-1} {\cal E}_{d_k} \right) \oplus {\cal E}_{d_{i}}^{\prime} \oplus \left( \bigoplus_{k=i+1}^{i-1} {\cal E}_{d_k} \right) \to \bigoplus_{k=1}^{l-1} p_{k,k+1} ^{*} {\cal O}_{CP^{N-1}} (N) \to 0, 
\end{equation}
where $p_{k,k+1}$ is a projection map defined on $\mathrm{Nor} (d_1 , \dots , d_l ; i)$:
\begin{equation}
p_{k, k+1} (\varphi_1 , \dots , \varphi_l ) := \varphi_k (\infty _k ) \; (= \varphi_{k+1} (0_{k+1} )).
\end{equation}

In the case of the limit part, we have to consider the part, $CP^{N-1} \times \overline{M}_{0,2|1} \simeq CP^{N-1}$, coming from the infinitessimal component $A$. On this component, we define a vector bundle ${\cal E}_{0}^{\prime}$ as ${\cal E}_{0}^{\prime} := {\cal O}_{CP^{N-1}} (N)$. Then we can extend the orbi-bundle ${\cal E}_{d}^{\prime}$ to $\mathrm{Lim} (d_1 , \dots , d_l ; i)$ by using the following exact sequence:
\begin{align}
0 &\to {\cal E}_{d}^{\prime}|_{\mathrm{Lim} (d_1 , \dots , d_l ; i)} \to \left( \bigoplus_{k=1}^{i} {\cal E}_{d_k} \right) \oplus {\cal E}_{0}^{\prime} \oplus \left( \bigoplus_{k=i+1}^{l} {\cal E}_{d_k} \right) \notag \\
&\to \left( \bigoplus_{k=1}^{i} p_{k,k+1} ^{*} {\cal O}_{CP^{N-1}} (N) \right) \oplus {\cal O}_{CP^{N-1}} (N) \oplus \left( \bigoplus_{k=i+1}^{l-1} p_{k,k+1} ^{*} {\cal O}_{CP^{N-1}} (N) \right) \to 0, \label{exseq}
\end{align}
where $p_{k,k+1}$ is a projection map defined on $\mathrm{Lim} (d_1 , \dots , d_l ; i)$:
\begin{equation}
p_{k, k+1} (\varphi_1 , \dots , \varphi_i , [\mathbf{a}] , \varphi_{i+1} , \dots ,  \varphi_l ) := \varphi_k (\infty_k ) =
\begin{cases}
\varphi_{k+1} (0_{k+1}) & (k=1, \dots , i-1, i+1 , \dots , l-1) \\
[\mathbf{a}] & (k=i).
\end{cases}
\end{equation}
In this way, we have extended rank $(Nd + 1)$ orbi-bundle ${\cal E}_{d}^{\prime}$ to the whole parts of $\widetilde{Mp}_{0,2|1}(N,d)$.

\section{Proof of the Main Theorem: Using Localization Technique}

In this section, we prove the main theorem. First, following \cite{JS}, we introduce the following $\mathbb{C}^{\times}$ action on $Mp_{0,2}(N,d)$:
\begin{equation}
e^t \cdot [ (\mathbf{a}_{0} , \dots , \mathbf{a}_{d} , z_1) ] = [(e^{\lambda_0 t} \mathbf{a}_{0} , \dots , e^{\lambda_d t} \mathbf{a}_{d} , z_1) ], \label{actloc}
\end{equation}
where $\lambda_0 , \dots , \lambda_d \in \mathbb{C}$ are the characters of $\mathbb{C}^{\times}$ action. 
We assume that $\lambda_{i}$'s are mutually different.
In order to extend the action to the whole parts of  $\widetilde{Mp}_{0,2|1}(N,d)$, Let us recall the case of $\widetilde{Mp}_{0,2}(N,d)$. In this case, $\mathbb{C}^{\times}$ action on a chain of quasimaps $(\varphi_{1},\cdots, \varphi_{l})\in Mp_{0,2} (N,d_1 ) \times_{CP^{N-1}} 
Mp_{0,2} (N,d_2 ) \times_{CP^{N-1}}\dotsb \times_{CP^{N-1}} Mp_{0,2} (N,d_l )$ is given as follows:
\ba
&&e^{t}\cdot([(\mathbf{a}_{0},\cdots,\mathbf{a}_{d_{1}})],\cdots,[(\mathbf{a}_{\sum_{j=1}^{l-1}d_{j}},\cdots,\mathbf{a}_{\sum_{j=1}^{l}d_{j}})])\no\\
&=&([(e^{\lambda_{0}t}\mathbf{a}_{0},\cdots,e^{\lambda_{d_{1}}t}\mathbf{a}_{d_{1}})],\cdots,[(e^{\lambda_{ \sum_{j=1}^{l-1}d_{j} }t}\mathbf{a}_{\sum_{j=1}^{l-1}d_{j}},\cdots,e^{\lambda_{ \sum_{j=1}^{l}d_{j} }t}\mathbf{a}_{\sum_{j=1}^{l}d_{j}})]).
\label{cext}
\ea
Then, a point in $\mathrm{Nor}(d_1 , \dots , d_l ; i)$ is represented by changing the $i$-th entry of the above expression into $[(\mathbf{a}_{\sum_{j=1}^{i-1}d_{j}},\cdots,\mathbf{a}_{\sum_{j=1}^{i}d_{j}},z_{1})]$. Hence we define the action by changing the above action on the $i$-th entry into,    
\ba
e^{t}\cdot[(\mathbf{a}_{\sum_{j=1}^{i-1}d_{j}},\cdots,\mathbf{a}_{\sum_{j=1}^{i}d_{j}},z_{1})]=[(e^{\lambda_{ \sum_{j=1}^{i-1}d_{j} }t}\mathbf{a}_{\sum_{j=1}^{i-1}d_{j}},\cdots,e^{\lambda_{ \sum_{j=1}^{i}d_{j} }t}\mathbf{a}_{\sum_{j=1}^{i}d_{j}},z_{1})].
\ea
On the other hand a point in $\mathrm{Lim}(d_1 , \dots , d_l ; i)$ is represented by inserting $([\mathbf{a}])$ between the $i$-th entry and the $(i+1)$-th entry of (\ref{cext}).
But this $CP^{N-1}$  facotor is trivialized by the fiber product and we don't need further consideration. 
     
We then determine the fixed point sets of  $\widetilde{Mp}_{0,2|1}(N,d)$ under the $\mathbb{C}^{\times}$ action.
By using the $\mathbb{C}^{\times}\times \mathbb{C}^{\times}$ action used in defining $Mp_{0,2|1}(N,d)$, the $\mathbb{C}^{\times}$ action (\ref{actloc}) is rewritten as,
\begin{align}
e^t \cdot [ (\mathbf{a}_{0} , \dots , \mathbf{a}_{i} , \dots , \mathbf{a}_{d} ,z_1 , \dots , z_n ] &= [e^{\lambda_0 t} \mathbf{a}_{0} , \dots , e^{\lambda_i t} \mathbf{a}_{i} , \dots , e^{\lambda_d t} \mathbf{a}_{d} , z_1) ] \notag \\
&= [(\mathbf{a}_{0} , \dots , e^{(\lambda_i - \lambda_0 ) t} \mathbf{a}_{i} , \dots , e^{( \lambda_d - \lambda_0 ) t} \mathbf{a}_{d} , z_1)] \notag \\
&= [(\mathbf{a}_{0} , \dots , e^{\left( \lambda_i - \lambda_0 - \frac{i (\lambda_d - \lambda_0 )t}{d} \right) t} \mathbf{a}_{i} , \dots , \mathbf{a}_{d} , e^{\frac{\lambda_d - \lambda_0 }{d} t}z_1) ].
\end{align}
Hence the fixed point set of $Mp_{0,2|1}(N,d)$ is empty because $e^{\frac{\lambda_d - \lambda_0 }{d} t}z_1\;\;(z_{1}\neq 0)$ is never fixed under variation of $t$.
But the fixed point set of $Mp_{0,2}(N,d)$ is not empty and it is given as follows. 
\begin{equation}
Z(d) := \left\{ [\mathbf{a}_0 , \dots , \mathbf{a}_{d}] \in Mp_{0,2}(N,d) \middle| \mathbf{a}_{1} = \dots = \mathbf{a}_{d-1} = \mathbf{0} \right\} .
\end{equation}
Note that $Z(d)$ is isomorphic to $CP^{N-1} \times CP^{N-1}$ and it is a set of $\mathbb{Z}_d$ singularities of $Mp_{0,2}(N,d)$. Furthermore, by the above discussion, fixed point sets of $\widetilde{Mp}_{0,2|1}(N,d)$ under the action (\ref{actloc}) lies the limit part $\mathrm{Lim}(d_1 , \dots , d_l ; i)$ and they are explicitly written as
\begin{align}
&Z(d_1 , \dots , d_i ; i) \notag \\ 
&:= Z(d_1 ) \times_{CP^{N-1}} \dots \times_{CP^{N-1}} Z(d_i ) \times_{CP^{N-1}} (CP^{N-1} \times \overline{M}_{0,2|1}) \times_{CP^{N-1}} Z(d_{i+1}) \times_{CP^{N-1}} \dots \times_{CP^{N-1}} Z(d_{l}) \notag \\
&\simeq (CP^{N-1})_{0} \times \dots \times (CP^{N-1})_{i} \times \overline{M}_{0,2|1} \times (CP^{N-1})_{i+1} \dots \times (CP^{N-1})_{l}.
\label{fixi}
\end{align}
At this stage, we introduce the following notation for brevity.
\ba
\Delta_{i}:=\sum_{j=1}^{i}d_{i}\;\;(\Delta_{0}:=0,\;\;\Delta_{l+1}=d (\Leftrightarrow d_{l+1}=0)).
\ea
Then $(CP^{N-1})_{k}\;(k=0,1,\cdots,l)$ in (\ref{fixi}) is image of the evaluation map,
\begin{equation}
\mbox{ev}_{0}^{(k)} : Mp_{0,2|0}(N,d_{k+1}) \to CP^{N-1} \; ; \; \left[(\mathbf{a}_{\Delta_{k}},\mathbf{a}_{\Delta_{k}+1},\cdots,\mathbf{a}_{\Delta_{k+1}})\right] \mapsto \mathbf{a}_{\Delta_{k}}.
\end{equation}
We note here that the locus $Z(d_1 , \dots , d_i ; i)$ are set of $\mathbb{Z}_{d_{1}}\times\cdots\times \mathbb{Z}_{d_{l}}$ orbifold singularities.

With these set-up's, we compute the intersection number $w \left( \sigma_{j} ({\cal O}_{h^{a}}) {\cal O}_{h^{b}} \middle| {\cal O}_{h} \right)_{0,d}$ by using localization technique.  
The resul is given as the following lemma. 
Recall that $\psi^{\prime}_{0}$ (resp. $\psi_{0}$) is the first Chern classs of the line bundle ${\cal L}^{\prime}_{0}$ on $\widetilde{Mp}_{0,2|1}(N,d)$ (resp. ${\cal L}_{0}$ on $\widetilde{Mp}_{0,2}(N,d)$) whose fiber at each point of the moduli space is given as cotangent space at the first marked point (${0}_{1}$) of the genus $0$ semi-stable curve. 
\begin{lem}
For $j > 0$, the following equality holds.
\begin{align}
&w \left( \sigma_{j} ({\cal O}_{h^{a}}) {\cal O}_{h^{b}} \middle| {\cal O}_{h} \right)_{0,d|1} \notag \\
&= \sum_{(d_1 , \dots , d_l ) \in S(d)} \frac{1}{d_1 \dots d_l} \sum_{k = 1}^{l} \int_{\prod _{i = 0}^{l} (CP^{N-1})_{i}} \Bigg\{ \bigl(\frac{h_{1}+\lambda_{\Delta_{1}}-h_{0}-\lambda_{\Delta_{0}}}{d_{1}}\bigr)^{j}   (h_0 + \lambda_{\Delta_0} )^a (h_l + \lambda_{\Delta_l})^b \notag \\
&\quad \cdot \frac{\prod_{i=1}^{l} \prod_{m_i = 0}^{Nd_i} \left( \frac{Nd_i - m_i}{d_i} (h_{i-1} + \lambda_{\Delta_{i-1}} ) + \frac{m_i}{d_i} (h_i + \lambda_{\Delta_{i}} ) \right)}{\prod_{i=1}^{l} \prod_{n_i = 1}^{d_i -1} \left( \frac{d_i - n_i}{d_i} (h_{i-1} + \lambda_{\Delta_{i-1}} ) + \frac{n_i}{d_i} ( h_i + \Delta_{i} ) - \lambda_{\Delta_{i-1} + n_i} \right)^{N}} \notag \\
&\quad \cdot \frac{\frac{d_k (h_k + \lambda_{\Delta_{k}})}{h_k + \lambda_{\Delta_{k}} - h_{k-1} - \lambda_{\Delta_{k-1}}} + \frac{d_{k+1} (h_k + \lambda_{\Delta_{k}})}{h_k + \lambda_{\Delta_k} - h_{k+1} - \lambda_{\Delta_{k+1}}}}{\prod_{i=1}^{l-1} \left( \frac{h_i + \lambda_{\Delta_i} - h_{i-1} - \lambda_{\Delta_{i-1}}}{d_i} + \frac{h_i + \lambda_{\Delta_{i}} - h_{i+1} - \lambda_{\Delta_{i+1}}}{d_{i+1}} \right) N(h_i + \lambda_{\Delta_{i}})} \Bigg\} , \label{lemma1}
\end{align}
where $\Delta_i := d_0 + \dots + d_i \ (d_0 = d_{l+1} := 0)$ and $h_i$ is the hyperplane class of $(CP^{N-1})_i$.
\end{lem}

\noindent
\textit{Proof.} \ We apply localization formula (for orbifold) to $w \left( \sigma_{j} ({\cal O}_{h^{a}}) {\cal O}_{h^{b}} \middle| {\cal O}_{h} \right)_{0,d}$:

\begin{align}
& \int_{\widetilde{Mp}_{0,2|1}(N,d)} (\psi_{0}^{\prime})^{j} \cdot \mbox{ev}_{0}^{*} (h^{a}) \cdot \mbox{ev}_{\infty}^{*} (h^{b}) \cdot \mbox{ev}_{1}^{*} (h) \cdot c_{top} ({\cal E}_{d}^{\prime}) \notag \\
&= \sum_{(d_1 , \dots , d_l ) \in S(d)} \sum_{k=0}^{l} \frac{1}{d_1 \dotsb d_l} \int_{Z(d_1 , \dots , d_l ; k)} \{ \bigl(c_1 ({\cal L}_{0}^{\prime} |_{Z(d_1 , \dots , d_l ; k)})^j  \cdot \tilde{c}_1 (\mbox{ev}_{0}^{*} (h)|_{Z(d_1 , \dots , d_l ; k)})^{a} \cdot \tilde{c}_1 (\mbox{ev}_{\infty}^{*} (h)|_{Z(d_1 , \dots , d_l ; k)})^{b} \notag \\
&\quad \cdot \tilde{c}_1 (\mbox{ev}_{1}^{*} (h)|_{Z(d_1 , \dots , d_l ; k)}) \cdot \tilde{c}_{Nd+1} ({\cal E}^{\prime}_{d}|_{Z(d_1 , \dots , d_l ; k)}) / \tilde{c}_{top} (N(d_1 , \dots , d_l ; k)) \} , \label{locint}
\end{align} 
where $h$ denotes hyperplane class of $CP^{N-1}$ and $N(d_1 , \dots , d_l ; k)$ denotes a normal bundle of the fixed point set $Z(d_1 , \dots , d_l ; k)$. $\tilde{c}_{i}(*)$ (resp. $\tilde{c}_{top}(*)$) means equivariant $i$-th (resp. top) Chern class of the vector bundle $*$. The factor $\frac{1}{d_1 \dotsb d_l}$ comes from dividing by cardinality of the group  $\mathbb{Z}_{d_{1}}\times\cdots\times \mathbb{Z}_{d_{l}}$.

First, we compute equivariant Chern classes in the right hand side. Recall that the domain of a quasimap in $\mathrm{Lim}(d_1 , \dots , d_l ; i)$ is a chain type semi-stable curve constructed by gluing $B_1,\cdots, B_{l}$ and $A$ (infinitessimal component) in the following order:
\begin{equation}
\begin{cases}
A \coprod \left(\coprod_{i=1}^{l} B_i \right) & (k=0) \\
\left( \coprod_{i=1}^{k} B_i \right)\coprod A \coprod\left(\coprod_{i=k+1}^{l} B_i \right)& (k = 1, \dots , l-1) \\
\left(  \coprod_{i=1}^{l} B_i \right) \coprod A & (k=l)
\end{cases}
\end{equation}
In the case of $B_{i}$, the situation is the same as the one of $Mp_{0,2}(N,d)$. So the compuation goes in the same way as the one given in \cite{Jin1, Jin3}.
 According to these literatures and \cite{JS}, the equivariant Chern classes $\tilde{c}_1 (\mbox{ev}_{0}^{*} (h)|_{Z(d_1 , \dots , d_l ; k)})$, $\tilde{c}_1 (\mbox{ev}_{\infty}^{*} (h)|_{Z(d_1 , \dots , d_l ; k)})$, $\tilde{c}_1 (\mbox{ev}_{1}^{*} (h)|_{Z(d_1 , \dots , d_l ; k)})$, and $\tilde{c}_{top} (N(d_1 , \dots , d_l ; k))$ are given by
\begin{equation}
\tilde{c}_1 (\mbox{ev}_{0}^{*} (h)|_{Z(d_1 , \dots , d_l ; k)}) = h_0 + \lambda_{\Delta_0},
\end{equation}
\begin{equation}
\tilde{c}_1 (\mbox{ev}_{\infty}^{*} (h)|_{Z(d_1 , \dots , d_l ; k)}) = h_l + \lambda_{\Delta_l},
\end{equation}
\begin{equation}
\tilde{c}_1 (\mbox{ev}_{1}^{*} (h)|_{Z(d_1 , \dots , d_l ; k)}) = h_k + \lambda_{\Delta_{k}},
\end{equation}
and
\begin{align}
\tilde{c}_{top} (N(d_1 , \dots , d_l ; k)) &= \left( \prod_{i=1}^{l} \prod_{n_1 = 1}^{d_1 - 1} \left( \frac{d_i - n_i }{d_i} (\lambda_{\Delta_{i-1}} + h_{i-1}) + \frac{n_i}{d_i} (\lambda_{\Delta_i} + h_i ) - \lambda_{\Delta_{i-1} + n_i} \right) ^N \right) \cdot s_{k}, \label{nor}
\end{align}
respectively. Here, $h_{i}$ denotes hyperplane class of $(CP^{N-1})_{i}$ in (\ref{fixi}). In (\ref{nor}), the first term 
\begin{equation}
\prod_{i=1}^{l} \prod_{n_1 = 1}^{d_1 - 1} \left( \frac{d_i - n_i }{d_i} (\lambda_{\Delta_{i-1}} + h_{i-1}) + \frac{n_i}{d_i} (\lambda_{\Delta_i} + h_i ) - \lambda_{\Delta_{i-1} + n_i} \right) ^N
\end{equation}
is the contribution from normal bundles of $Z(d_1 ), \dots , Z(d_l )$.  On the other hand, the foloowing factors: 
\begin{equation}
s_0 := 
\left( \frac{h_0 + \lambda_{\Delta_0} - h_1 - \lambda_{\Delta_1}}{d_1} + c_1 ( T_{\infty^{\prime}}^{\prime} A) \right) \left( \prod_{i = 1}^{l-1} \left( \frac{h_i + \lambda_{\Delta_i} - h_{k-1} - \lambda_{\Delta_{i-1}}}{d_i} + \frac{h_i + \lambda_{\Delta_i} - h_{i+1} - \lambda_{\Delta_{i+1}}}{d_{i+1}} \right) \right),
\end{equation}
\begin{align}
s_k &:= \left( \prod_{i = 1}^{k-1} \left( \frac{h_i + \lambda_{\Delta_i} - h_{k-1} - \lambda_{\Delta_{i-1}}}{d_i} + \frac{h_i + \lambda_{\Delta_i} - h_{i+1} - \lambda_{\Delta_{i+1}}}{d_{i+1}} \right) \right) \notag \\
&\ \quad \cdot \left( \frac{h_k + \lambda_{\Delta_k} - h_{k-1} - \lambda_{\Delta_{k-1}}}{d_k} + c_1 ( T_{0^{\prime}}^{\prime} A) \right) \left( \frac{h_k + \lambda_{\Delta_k} - h_{k+1} - \lambda_{\Delta_{k+1}}}{d_{k+1}} + c_1 ( T_{\infty^{\prime}}^{\prime} A) \right) \notag \\
&\ \quad \cdot \left( \prod_{i = k+1}^{l-1} \left( \frac{h_i + \lambda_{\Delta_i} - h_{k-1} - \lambda_{\Delta_{i-1}}}{d_i} + \frac{h_i + \lambda_{\Delta_i} - h_{i+1} - \lambda_{\Delta_{i+1}}}{d_{i+1}} \right) \right) \quad (k = 1, \dots , l-1),
\end{align}
\begin{equation}
s_l := \left( \prod_{i = 1}^{l-1} \left( \frac{h_i + \lambda_{\Delta_i} - h_{k-1} - \lambda_{\Delta_{i-1}}}{d_i} + \frac{h_i + \lambda_{\Delta_i} - h_{i+1} - \lambda_{\Delta_{i+1}}}{d_{i+1}} \right) \right) \left( \frac{h_l + \lambda_{\Delta_l} - h_{l-1} - \lambda_{\Delta_{l-1}}}{d_l} + c_1 ( T_{0^{\prime}}^{\prime} A) \right)
\end{equation}
come from smoothing nodal singularities (see \cite{Jin1, JS}). Next, we compute the remaining equivalent Chern classes i.e., $\tilde{c}_1 ({\cal L}_{0}^{\prime} |_{Z(d_1 , \dots , d_l ; k)})$ and $\tilde{c}_{Nd+1} ({\cal E}_{d}^{\prime}|_{Z(d_1 , \dots , d_l ; k)})$. Since a fiber of the line bundle ${\cal L}_{0}^{\prime} |_{Z(d_1 , \dots , d_l ; k)}$ at a point in $Z(d_1 , \dots , d_l ; k)$ doesn't depend on coefficient vectors of quasimaps, the first equivarlant Chern class reduces to ordinary first Chern class $c_1 ({\cal L}_{0}^{\prime} |_{Z(d_1 , \dots , d_l ; k)})$. 
On the other hand, ${\cal E}_{d}^{\prime}|_{Z(d_1 , \dots , d_l ; k)}$ doesn't depend on the informaion of $z_{1}$ at all. Therefore, we can use the discussion given in \cite{Jin1, Jin3}
(one subtlety comes from the contribution ${\cal E}_{0}^{\prime}$ in (\ref{exseq}), but it will soon turn out to be irrelevant). 
By using the exact sequence (\ref{exseq}) and the discussions in \cite{Jin1,Jin3}, $\tilde{c}_{Nd+1} ({\cal E}_{d}^{\prime}|_{Z(d_1 , \dots , d_l ; k)})$ is explicitly given as follows:
\begin{align}
\tilde{c}_{Nd+1} ({\cal E}_{d}^{\prime}|_{Z(d_1 , \dots , d_l ; k)}) &= \Bigg\{ \left( \prod_{i=1}^{k} \prod_{m_i = 0}^{Nd_i} \left( \frac{(Nd_i - m_i )(h_{i-1} + \lambda_{\Delta_{i-1}}) + m_i (h_i + \lambda_{\Delta_{i}})}{d_i} \right) \right) \cdot \tilde{c}_1 ({\cal E}_0) \notag \\
&\qquad \cdot \left( \prod_{i=k+1}^{l} \prod_{m_i = 0}^{Nd_i} \left( \frac{(Nd_i - m_i )(h_{i-1} + \lambda_{\Delta_{i-1}}) + m_i (h_i + \lambda_{\Delta_{i}})}{d_i} \right) \right) \Bigg\} \notag \\
&\ \ \quad \Bigg/ \left\{ \left( \prod_{i=1}^{k} N(h_i + \lambda_{\Delta_{i}}) \right) \cdot \tilde{c}_1 ({\cal O}_{CP^{N-1}} (N)) \cdot \left( \prod_{i=k+1}^{l-1} N(h_i + \lambda_{\Delta_{i}}) \right) \right\} \notag \\
&= \frac{\prod_{i=1}^{l} \prod_{m_i = 0}^{Nd_i} \left( \frac{(Nd_i - m_i )(h_{i-1} + \lambda_{\Delta_{i-1}}) + m_i (h_i + \lambda_{\Delta_{i}})}{d_i} \right)}{\prod_{i=1}^{l-1} N(h_i + \lambda_{\Delta_{i}})}.
\end{align}
In going from the r.h.s. of the first line to the second line, we used ${\cal E}_0 = {\cal O}_{CP^{N-1}} (N)$. Therefore, the integrand in the right hand side of (\ref{locint}) is given by,
\begin{equation}
\frac{c_1 ({\cal L}_{0}^{\prime} |_{Z(d_1 , \dots , d_l ; k)})^j  (h_0 + \lambda_{\Delta_0})^a (h_l + \lambda_{\Delta_l})^b (h_k + \lambda_{\Delta_{k}}) \frac{\prod_{i=1}^{l} \prod_{m_i = 0}^{Nd_i} \left( \frac{(Nd_i - m_i )(h_{i-1} + \lambda_{\Delta_{i-1}}) + m_i (h_i + \lambda_{\Delta_{i}})}{d_i} \right)}{\prod_{i=1}^{l-1} N(h_i + \lambda_{\Delta_{i}})}}{\left( \prod_{i=1}^{l} \prod_{n_1 = 1}^{d_1 - 1} \left( \frac{d_i - n_i }{d_i} (\lambda_{\Delta_{i-1}} + h_{i-1}) + \frac{n_i}{d_i} (\lambda_{\Delta_i} + h_i ) - \lambda_{\Delta_i + n_i} \right) ^N \right) s_{k}}.
\label{integ}
\end{equation}
At this stage, let us discuss the  part $c_1 ({\cal L}_{0}^{\prime} |_{Z(d_1 , \dots , d_l ; k)})^j $, Recall that fiber of ${\cal L}_{0}^{\prime}$ is given as cotangent space at $0_{1}$ of the first component of $CP^{1}$ if $k\geq 1$ and from cotangent space of $0^{\prime}$ of $A$ if $k=0$. Hence we obtain, 
\ba
c_1 ({\cal L}_{0}^{\prime} |_{Z(d_1 , \dots , d_l ; k)})^j =\left\{ \begin{array}{ll}c_1 (T_{0^{\prime}}^{\prime *} A)^{j}& (k=0)\\
&\\
 \bigl(\frac{h_{1}+\lambda_{\Delta_{1}}-h_{0}-\lambda_{\Delta_{0}}}{d_{1}}\bigr)^{j} &(k=1,2,\cdots,l),    \end{array}\right.                                                                                               
\ea
by using the discussion given in \cite{JS}.  
At this stage, the term that depends on $A$ in (\ref{integ}) is given as follows.
\begin{equation}
I_A (k) :=
\begin{cases}
c_1 (T_{0^{\prime}}^{\prime *} A)^{j}\bigg/ \left( \frac{h_0 + \lambda_{\Delta_0} - h_1 - \lambda_{\Delta_1}}{d_1} + c_1 ( T_{\infty^{\prime}}^{\prime} A) \right) \quad (k = 0) \\
1 \bigg/ \left\{ \left( \frac{h_k + \lambda_{\Delta_k} - h_{k-1} - \lambda_{\Delta_{k-1}}}{d_k} + c_1 ( T_{0^{\prime}}^{\prime} A) \right) \left( \frac{h_k + \lambda_{\Delta_k} - h_{k+1} - \lambda_{\Delta_{k+1}}}{d_{k+1}} + c_1 ( T_{\infty^{\prime}}^{\prime} A) \right) \right\} \\
\quad (k=1 , \dots , l-1) \\
1 \bigg/ \left( \frac{h_l + \lambda_{\Delta_l} - h_{l-1} - \lambda_{\Delta_{l-1}}}{d_l} + c_1 ( T_{0^{\prime}}^{\prime} A) \right) \quad (k = l). 
\end{cases}
\end{equation}
To integrate $I_A (k)$ on $\overline{M}_{0,2|1}$, one can use the following formula \cite{AG, MOP}:
\begin{equation}
\int_{\overline{M}_{0,2|1}} (c_1 (T_{0^{\prime}}^{\prime *} A))^p (c_1 (T_{\infty^{\prime}}^{\prime *} A))^q = 
\begin{cases}
1 & (p=q=0) \\
0 & (otherwise).
\end{cases}
\end{equation}
Alternatively, the above formula can be derived from the fact that $\overline{M}_{0,2|1}$ is just a point.
 Hence integration of $I_A (0)$ on $\overline{M}_{0,2|1}$ goes as follows.
\begin{align}
\int_{\overline{M}_{0,2|1}} I_A (0) &= \int_{\overline{M}_{0,2|1}} \frac{c_1 (T_{0^{\prime}}^{\prime *} A)^j}{\frac{h_0 + \lambda_{\Delta_0} - h_1 - \lambda_{\Delta_1}}{d_1} + c_1 ( T_{\infty^{\prime}}^{\prime} A)} \notag \\
&= \sum_{q=0}^{\infty} \frac{1}{\left( \frac{h_0 + \lambda_{\Delta_0} - h_1 - \lambda_{\Delta_1}}{d_1} \right)^{q+1}} \int_{\overline{M}_{0,2|1}} c_1 (T_{0^{\prime}}^{\prime *} A)^j (c_1 (T_{\infty^{\prime}}^{\prime *} A))^q \notag \\
&=0,
\end{align}
where we used the fact that $j$ is positive. In the cases of $k=1, \dots , l-1$, we obtain,
\begin{align}
\int_{\overline{M}_{0,2|1}} I_A (k) &= \int_{\overline{M}_{0,2|1}} \frac{1}{\left( \frac{h_k + \lambda_{\Delta_k} - h_{k-1} - \lambda_{\Delta_{k-1}}}{d_k} + c_1 ( T_{0^{\prime}}^{\prime} A) \right) \left( \frac{h_k + \lambda_{\Delta_k} - h_{k+1} - \lambda_{\Delta_{k+1}}}{d_{k+1}} + c_1 ( T_{\infty^{\prime}}^{\prime} A) \right)} \notag \\
&= \sum_{p, q=0}^{\infty} \frac{1}{\left( \frac{h_k + \lambda_{\Delta_k} - h_{k-1} - \lambda_{\Delta_{k-1}}}{d_k} \right)^{p+1} \left( \frac{h_k + \lambda_{\Delta_k} - h_{k+1} - \lambda_{\Delta_{k+1}}}{d_{k+1}} \right)^{q+1}} \int_{\overline{M}_{0,2|1}} (c_1 (T_{0^{\prime}}^{\prime *} A))^p (c_1 (T_{\infty^{\prime}}^{\prime *} A))^q \notag \\
&=  \frac{d_k}{h_k + \lambda_{\Delta_k} - h_{k-1} - \lambda_{\Delta_{k-1}}} \cdot \frac{d_{k+1}}{h_k + \lambda_{\Delta_k} - h_{k+1} - \lambda_{\Delta_{k+1}}} \cdot 1 \notag \\
&= \frac{\frac{d_k}{h_k + \lambda_{\Delta_k} - h_{k-1} - \lambda_{\Delta_{k-1}}} + \frac{d_{k+1}}{h_k + \lambda_{\Delta_k} - h_{k+1} - \lambda_{\Delta_{k+1}}}}{\frac{h_k + \lambda_{\Delta_k} - h_{k-1} - \lambda_{\Delta_{k-1}}}{d_k} + \frac{h_k + \lambda_{\Delta_k} - h_{k+1} - \lambda_{\Delta_{k+1}}}{d_{k+1}}}.
\end{align}
Integration in the case of $k=l$ is given by,
\begin{align}
\int_{\overline{M}_{0,2|1}} I_A (l) &= \int_{\overline{M}_{0,2|1}} \frac{1}{\frac{h_l + \lambda_{\Delta_l} - h_{l-1} - \lambda_{\Delta_{l-1}}}{d_l} + c_1 ( T_{0^{\prime}}^{\prime} A)} \notag \\
&= \sum_{p=0}^{\infty} \frac{1}{\left( \frac{h_l + \lambda_{\Delta_l} - h_{l-1} - \lambda_{\Delta_{l-1}}}{d_l} \right)^{p+1}} \int_{\overline{M}_{0,2|1}} c_1 (T_{0^{\prime}}^{\prime *} A)^p \notag 
\\
&= \frac{d_l}{h_l + \lambda_{\Delta_l} - h_{l-1} - \lambda_{\Delta_{l-1}}}.
\end{align}
At this stage, the remaining domain of integration is the product space $\prod_{i=1}^{l} (CP^{N-1})_i$. Therefore, we obatin the 
following equalities:
\begin{align}
&w \left( \sigma_{j} ({\cal O}_{h^{a}}) {\cal O}_{h^{b}} \middle| {\cal O}_{h} \right)_{0,d} \notag \\
&= \sum_{(d_1 , \dots , d_l ) \in S(d)} \frac{1}{d_1 \dots d_l} \Bigg\{ \notag \\
&\quad 0 + \sum_{k=1}^{l-1} \int_{\prod_{i=1}^{l} (CP^{N-1})_i} \bigl(\frac{h_{1}+\lambda_{\Delta_{1}}-h_{0}-\lambda_{\Delta_{0}}}{d_{1}}\bigr)^{j} (h_0 + \lambda_{\Delta_0})^a (h_l + \lambda_{\Delta_l})^b (h_k + \lambda_{\Delta_{k}}) \notag \\
&\qquad \qquad \cdot \frac{ \frac{\prod_{i=1}^{l} \prod_{m_i = 0}^{Nd_i} \left( \frac{(Nd_i - m_i )(h_{i-1} + \lambda_{\Delta_{i-1}}) + m_i (h_i + \lambda_{\Delta_{i}})}{d_i} \right)}{\prod_{i=1}^{l-1} N(h_i + \lambda_{\Delta_{i}})}}{\left( \prod_{i=1}^{l} \prod_{n_1 = 1}^{d_1 - 1} \left( \frac{d_i - n_i }{d_i} (\lambda_{\Delta_{i-1}} + h_{i-1}) + \frac{n_i}{d_i} (\lambda_{\Delta_i} + h_i ) - \lambda_{\Delta_i + n_i} \right) ^N \right)} \notag \\
&\qquad \qquad \cdot \frac{\frac{d_k}{h_k + \lambda_{\Delta_k} - h_{k-1} - \lambda_{\Delta_{k-1}}} + \frac{d_{k+1}}{h_k + \lambda_{\Delta_k} - h_{k+1} - \lambda_{\Delta_{k+1}}}}{ \prod_{i=1}^{l-1} \left( \frac{h_i + \lambda_{\Delta_i} - h_{i-1} - \lambda_{\Delta_{i-1}}}{d_i} + \frac{h_i + \lambda_{\Delta_i} - h_{i+1} - \lambda_{\Delta_{i+1}}}{d_{i+1}} \right)} \notag \\
&\quad + \int_{\prod_{i=1}^{l} (CP^{N-1})_i}  \bigl(\frac{h_{1}+\lambda_{\Delta_{1}}-h_{0}-\lambda_{\Delta_{0}}}{d_{1}}\bigr)^{j} (h_0 + \lambda_{\Delta_0})^a (h_l + \lambda_{\Delta_l})^b (h_l + \lambda_{\Delta_{l}}) \notag \\
&\qquad \cdot \frac{ \frac{\prod_{i=1}^{l} \prod_{m_i = 0}^{Nd_i} \left( \frac{(Nd_i - m_i )(h_{i-1} + \lambda_{\Delta_{i-1}}) + m_i (h_i + \lambda_{\Delta_{i}})}{d_i} \right)}{\prod_{i=1}^{l-1} N(h_i + \lambda_{\Delta_{i}})}}{\left( \prod_{i=1}^{l} \prod_{n_1 = 1}^{d_1 - 1} \left( \frac{d_i - n_i }{d_i} (\lambda_{\Delta_{i-1}} + h_{i-1}) + \frac{n_i}{d_i} (\lambda_{\Delta_i} + h_i ) - \lambda_{\Delta_i + n_i} \right) ^N \right)} \notag \\
&\qquad \cdot \frac{\frac{d_l}{h_l + \lambda_{\Delta_l} - h_{l-1} - \lambda_{\Delta_{l-1}}}}{\prod_{i=1}^{l-1} \left( \frac{h_i + \lambda_{\Delta_i} - h_{i-1} - \lambda_{\Delta_{i-1}}}{d_i} + \frac{h_i + \lambda_{\Delta_i} - h_{i+1} - \lambda_{\Delta_{i+1}}}{d_{i+1}} \right)} \Bigg\} \notag \\
&= \sum_{(d_1 , \dots , d_l ) \in S(d)} \frac{1}{d_1 \dots d_l} \sum_{k=1}^{l} \int_{\prod_{i=1}^{l} (CP^{N-1})_i} \bigl(\frac{h_{1}+\lambda_{\Delta_{1}}-h_{0}-\lambda_{\Delta_{0}}}{d_{1}}\bigr)^{j} (h_0 + \lambda_{\Delta_0})^a (h_l + \lambda_{\Delta_l})^b \notag \\
&\qquad \qquad \cdot \frac{\prod_{i=1}^{l} \prod_{m_i = 0}^{Nd_i} \left( \frac{(Nd_i - m_i )(h_{i-1} + \lambda_{\Delta_{i-1}}) + m_i (h_i + \lambda_{\Delta_{i}})}{d_i} \right)}{\prod_{i=1}^{l} \prod_{n_1 = 1}^{d_1 - 1} \left( \frac{d_i - n_i }{d_i} (\lambda_{\Delta_{i-1}} + h_{i-1}) + \frac{n_i}{d_i} (\lambda_{\Delta_i} + h_i ) - \lambda_{\Delta_i + n_i} \right) ^N} \notag \\
&\qquad \qquad \cdot \frac{\frac{d_k (h_k + \lambda_{\Delta_{k}})}{h_k + \lambda_{\Delta_k} - h_{k-1} - \lambda_{\Delta_{k-1}}} + \frac{d_{k+1} (h_k + \lambda_{\Delta_{k}})}{h_k + \lambda_{\Delta_k} - h_{k+1} - \lambda_{\Delta_{k+1}}}}{ \prod_{i=1}^{l-1} \left( \frac{h_i + \lambda_{\Delta_i} - h_{i-1} - \lambda_{\Delta_{i-1}}}{d_i} + \frac{h_i + \lambda_{\Delta_i} - h_{i+1} - \lambda_{\Delta_{i+1}}}{d_{i+1}} \right) N(h_i + \lambda_{\Delta_i})},
\end{align}
and the last expression is nothing but the right hand side of (\ref{lemma1}). $\Box$

\vspace{1em}

Next, we state the following lemma:
\begin{lem}
\begin{equation}
\sum_{k=1}^{l} \left( \frac{d_k (h_k + \lambda_{\Delta_{k}})}{h_k + \lambda_{\Delta_{k}} - h_{k-1} - \lambda_{\Delta_{k-1}}} + \frac{d_{k+1} (h_k + \lambda_{\Delta_{k}})}{h_k + \lambda_{\Delta_{k}} - h_{k+1} - \lambda_{\Delta_{k+1}}} \right) = d + \frac{d_1 (h_0 + \lambda_{\Delta_{0}})}{h_1 + \lambda_{\Delta_{1}} - h_0 - \lambda_{\Delta_{0}}}.
\end{equation}
\end{lem}
\textit{Proof.} \ The proof of this lemma is straightforward. \ $\square$

\vspace{2em}

\noindent 
({\bf Proof of the Main Theorem})\quad     Combining the above two lemmas, $w \left( \sigma_{j} ({\cal O}_{h^{a}}) {\cal O}_{h^{b}} \middle| {\cal O}_{h} \right)_{0,d}$ is rewritten as,
\begin{align}
w \left( \sigma_{j} ({\cal O}_{h^{a}}) {\cal O}_{h^{b}} \middle| {\cal O}_{h} \right)_{0,d} 
&= d \sum_{(d_1 , \dots , d_l ) \in S(d)} \frac{1}{d_1 \dots d_l} \int_{\prod _{i = 0}^{l} (CP^{N-1})_{i}} \Bigg\{ \bigl(\frac{h_{1}+\lambda_{\Delta_{1}}-h_{0}-\lambda_{\Delta_{0}}}{d_{1}}\bigr)^{j} (h_0 + \lambda_{\Delta_0} )^a (h_l + \lambda_{\Delta_l})^b \notag \\
&\quad \cdot \frac{\prod_{i=1}^{l} \prod_{m_i = 0}^{Nd_i} \left( \frac{Nd_i - m_i}{d_i} (h_{i-1} + \lambda_{\Delta_{i-1}} ) + \frac{m_i}{d_i} (h_i + \lambda_{\Delta_{i}} ) \right)}{\prod_{i=1}^{l} \prod_{n_i = 1}^{d_i -1} \left( \frac{d_i - n_i}{d_i} (h_{i-1} + \lambda_{\Delta_{i-1}} ) + \frac{n_i}{d_i} ( h_i + \Delta_{i} ) - \lambda_{\Delta_{i-1} + n_i} \right)^{N}} \notag \\
&\quad \cdot \frac{1}{\prod_{i=1}^{l-1} \left( \frac{h_i + \lambda_{\Delta_i} - h_{i-1} - \lambda_{\Delta_{i-1}}}{d_i} + \frac{h_i + \lambda_{\Delta_{i}} - h_{i+1} - \lambda_{\Delta_{i+1}}}{d_{i+1}} \right) N(h_i + \lambda_{\Delta_{i}})} \Bigg\} \notag \\ 
&\quad + \sum_{(d_1 , \dots , d_l ) \in S(d)} \frac{1}{d_1 \dots d_l} \int_{\prod _{i = 0}^{l} (CP^{N-1})_{i}} \Bigg\{ \bigl(\frac{h_{1}+\lambda_{\Delta_{1}}-h_{0}-\lambda_{\Delta_{0}}}{d_{1}}\bigr)^{j} (h_0 + \lambda_{\Delta_0} )^a (h_l + \lambda_{\Delta_l})^b \notag \\
&\quad \cdot \frac{\prod_{i=1}^{l} \prod_{m_i = 0}^{Nd_i} \left( \frac{Nd_i - m_i}{d_i} (h_{i-1} + \lambda_{\Delta_{i-1}} ) + \frac{m_i}{d_i} (h_i + \lambda_{\Delta_{i}} ) \right)}{\prod_{i=1}^{l} \prod_{n_i = 1}^{d_i -1} \left( \frac{d_i - n_i}{d_i} (h_{i-1} + \lambda_{\Delta_{i-1}} ) + \frac{n_i}{d_i} ( h_i + \Delta_{i} ) - \lambda_{\Delta_{i-1} + n_i} \right)^{N}} \notag \\
&\quad \cdot \frac{1}{\prod_{i=1}^{l-1} \left( \frac{h_i + \lambda_{\Delta_i} - h_{i-1} - \lambda_{\Delta_{i-1}}}{d_i} + \frac{h_i + \lambda_{\Delta_{i}} - h_{i+1} - \lambda_{\Delta_{i+1}}}{d_{i+1}} \right) N(h_i + \lambda_{\Delta_{i}})} \notag \\
&\quad \cdot \frac{d_1 (h_0 + \lambda_{\Delta_{0}})}{h_1 + \lambda_{\Delta_{1}} - h_0 - \lambda_{\Delta_{0}}} \Bigg\}. \label{conc1}
\end{align}

On the other hand, we can compute $w \left( \sigma_{j} ({\cal O}_{h^{a}}) {\cal O}_{h^{b}} \right)_{0,d}$ by applyong the same consideration so far as follows,
\begin{align}
w \left( \sigma_{j} ({\cal O}_{h^{a}}) {\cal O}_{h^{b}} \right)_{0,d} &= \sum_{(d_1 , \dots , d_l ) \in S(d)} \frac{1}{d_1 \dots d_l} \sum_{k = 1}^{l} \int_{\prod _{i = 0}^{l} (CP^{N-1})_{i}} \Bigg\{ \bigl(\frac{h_{1}+\lambda_{\Delta_{1}}-h_{0}-\lambda_{\Delta_{0}}}{d_{1}}\bigr)^{j} (h_0 + \lambda_{\Delta_0} )^a (h_l + \lambda_{\Delta_l})^b \notag \\
&\quad \cdot \frac{\prod_{i=1}^{l} \prod_{m_i = 0}^{Nd_i} \left( \frac{Nd_i - m_i}{d_i} (h_{i-1} + \lambda_{\Delta_{i-1}} ) + \frac{m_i}{d_i} (h_i + \lambda_{\Delta_{i}} ) \right)}{\prod_{i=1}^{l} \prod_{n_i = 1}^{d_i -1} \left( \frac{d_i - n_i}{d_i} (h_{i-1} + \lambda_{\Delta_{i-1}} ) + \frac{n_i}{d_i} ( h_i + \Delta_{i} ) - \lambda_{\Delta_{i-1} + n_i} \right)^{N}} \notag \\
&\quad \cdot \frac{1}{\prod_{i=1}^{l-1} \left( \frac{h_i + \lambda_{\Delta_i} - h_{i-1} - \lambda_{\Delta_{i-1}}}{d_i} + \frac{h_i + \lambda_{\Delta_{i}} - h_{i+1} - \lambda_{\Delta_{i+1}}}{d_{i+1}} \right) N(h_i + \lambda_{\Delta_{i}})} \Bigg\} . \label{2ndterm}
\end{align}
Of course, this formula gives the same intersection number as the formula used in \cite{JM}:
\ba
w \left( \sigma_{j} ({\cal O}_{h^{a}}) {\cal O}_{h^{b}} \right)_{0,d} &=
&\frac{1}{(2 \pi \sqrt{-1})^{d+1}} \oint_{C_{0}} \frac{dz_{0}}{(z_{0})^{N}} \oint_{C_{1}} \frac{dz_{1}}{(z_{1})^{N}} \dots \oint_{C_{d}} \frac{dz_{d}}{(z_{d})^{N}} \no\\
&& \times (z_{0})^{a} (z_{1} - z_{0})^{j} \prod_{l=1}^{d}\bigl(\prod_{j=1}^{N}(jz_{l-1}+(N-j)z_{l})\bigr) \prod_{l=1}^{d-1} \frac{1}{Nz_{l} (2z_{l} - z_{l-1} - z_{l+1})} \no\\
&& \times (z_{d})^{b}.
\label{residue}
\ea
We can derive (\ref{residue}) by rewriting (\ref{2ndterm}) into residue integral by using the correpondence:
\ba
\int_{(CP^{N-1})}(h_{i})^{a}\;\leftrightarrow\;\frac{1}{2\pi\sqrt{-1}}\oint_{C_{0}}\frac{dz_{i}}{(z_{i})^{N}}(z_{i})^a,
\ea
and applying the standard procedure including non-equivariant limit ($\lambda_{*}\rightarrow  0$), 
which is presented in \cite{Jin1,Jin3}.

Hence the first term in the right hand side of (\ref{conc1}) is just $d \cdot w \left( \sigma_{j} ({\cal O}_{h^{a}}) {\cal O}_{h^{b}} \right)_{0,d}$. As for the second term of (\ref{conc1}), we have, 
\begin{align} 
(\text{the second term}) &= \sum_{(d_1 , \dots , d_l ) \in S(d)} \frac{1}{d_1 \dots d_l} \int_{\prod _{i = 0}^{l} (CP^{N-1})_{i}} \Bigg\{ \bigl(\frac{h_{1}+\lambda_{\Delta_{1}}-h_{0}-\lambda_{\Delta_{0}}}{d_{1}}\bigr)^{j-1} (h_0 + \lambda_{\Delta_0} )^{a+1}a (h_l + \lambda_{\Delta_l})^b \notag \\
&\quad \cdot \frac{\prod_{i=1}^{l} \prod_{m_i = 0}^{Nd_i} \left( \frac{Nd_i - m_i}{d_i} (h_{i-1} + \lambda_{\Delta_{i-1}} ) + \frac{m_i}{d_i} (h_i + \lambda_{\Delta_{i}} ) \right)}{\prod_{i=1}^{l} \prod_{n_i = 1}^{d_i -1} \left( \frac{d_i - n_i}{d_i} (h_{i-1} + \lambda_{\Delta_{i-1}} ) + \frac{n_i}{d_i} ( h_i + \Delta_{i} ) - \lambda_{\Delta_{i-1} + n_i} \right)^{N}} \notag \\
&\quad \cdot \frac{1}{\prod_{i=1}^{l-1} \left( \frac{h_i + \lambda_{\Delta_i} - h_{i-1} - \lambda_{\Delta_{i-1}}}{d_i} + \frac{h_i + \lambda_{\Delta_{i}} - h_{i+1} - \lambda_{\Delta_{i+1}}}{d_{i+1}} \right) N(h_i + \lambda_{\Delta_{i}})}\Bigg\},
\end{align}
by using,
\ba
\bigl(\frac{h_{1}+\lambda_{\Delta_{1}}-h_{0}-\lambda_{\Delta_{0}}}{d_{1}}\bigr)^{j} \cdot \frac{d_1 (h_0 + \lambda_{\Delta_{0}})}{h_1 + \lambda_{\Delta_{1}} - h_0 - \lambda_{\Delta_{0}}}
=\bigl(\frac{h_{1}+\lambda_{\Delta_{1}}-h_{0}-\lambda_{\Delta_{0}}}{d_{1}}\bigr)^{j-1}\cdot (h_0 + \lambda_{\Delta_{0}}).
\ea 
Looking back at (\ref{2ndterm}), we can see that the second term coincides with the intersection number $w \left( \sigma_{j-1} ({\cal O}_{h^{a+1}}) {\cal O}_{h^{b}} \right)_{0,d}$. 
Hence we have completed the proof of main theorem. $\Box$

\begin{Rem}
Let ${\cal E}_{d,\infty}^{\prime}$ (resp. ${\cal E}_{d,\infty}$) be orbi-bundle on $\widetilde{Mp}_{0,2|1}(N,d)$ (resp. $\widetilde{Mp}_{0,2}(N,d)$) defined by the following 
exact sequences \cite{CK, zinger}.
\ba
&&0\longrightarrow {\cal E}_{d,\infty}^{\prime}\longrightarrow {\cal E}_{d}^{\prime}\longrightarrow \mbox{ev}^{*}_{\infty}{\cal O}_{CP^{N-1}}(N)\longrightarrow  0,\no\\
&&0\longrightarrow {\cal E}_{d,\infty}\longrightarrow {\cal E}_{d}\longrightarrow \mbox{ev}^{*}_{\infty}{\cal O}_{CP^{N-1}}(N)\longrightarrow  0.
\ea
Then the intersection numbers $w \left( \sigma_{j} ({\cal O}_{h^{N-2-j}}) {\cal O}_{h^{-1}} \middle| {\cal O}_{h} \right)_{0,d}$ and 
$w \left( \sigma_{j} ({\cal O}_{h^{N-2-j}}) {\cal O}_{h^{-1}} \right)_{0,d}$ are defind as follows.  
\begin{align}
&\frac{1}{N}w \left( \sigma_{j} ({\cal O}_{h^{N-2-j}}) {\cal O}_{h^{-1}} \middle| {\cal O}_{h} \right)_{0,d} \notag \\
&:= \int_{\widetilde{Mp}_{0,2|1}(N,d)} (\psi^{\prime}_{0})^{j} \cdot \mbox{ev}_{0}^{*} (h^{N-2-j}) \cdot\mbox{ev}_{\infty}^{*} (h^{0})\cdot \mbox{ev}_{1}^{*} (h) \cdot c_{top} ({\cal E}_{d,\infty}^{\prime}),\notag\\
&\frac{1}{N}w \left( \sigma_{j} ({\cal O}_{h^{N-2-j}}) {\cal O}_{h^{-1}} \right)_{0,d} \notag \\
&:= \int_{\widetilde{Mp}_{0,2}(N,d)} (\psi_{0})^{j} \cdot \mbox{ev}_{0}^{*} (h^{N-2-j}) \cdot \mbox{ev}_{\infty}^{*} (h^{0})\cdot c_{top} ({\cal E}_{d,\infty}).
\end{align}
Then proof of the main theorem in the case of $b=-1$ is given by replacing ${\cal E}_{d}^{\prime}$ and ${\cal E}_{d}$ in the discussions so far by 
${\cal E}_{d,\infty}^{\prime}$ and ${\cal E}_{d,\infty}$. 
\label{final}  
\end{Rem}

\begin{Rem}
By changing orbi-bundles ${\cal E}_{d}^{\prime}$ and ${\cal E}_{d}$ into the ones obtained from degree $k$ homogeneous polynomial of $CP^{N-1}$, we can straightforwardly extend 
the main theorem to the case of degree $k$ hypersurface in $CP^{N-1}$. 
\end{Rem}


\begin{thebibliography}{99}
\bibitem{AG} V. Alexeev, G. Michael Guy, \textit{Moduli of weighted stable maps and their gravitational descendants}, J. Inst. Math. Jussieu \textbf{7}(3) (2008), 425--456.
\bibitem{CJ} A. Collino, M.Jinzenji.  \textit{On the structure of the small quantum cohomology rings of projective hypersurfaces.} Comm. Math. Phys. 206 (1999), no. 1, 157--183.  
\bibitem{CK} D.-A. Cox, S. Katz, \textit{ Mirror symmetry and algebraic geometry. }  Mathematical Surveys and Monographs, 68. American Mathematical Society, Providence, RI, 1999. 469 pp. ISBN: 0--8218--1059--6.
\bibitem{giv} A. B. Givental. \textit{Equivariant Gromov-Witten invariants.} Internat. Math. Res. Notices 1996, no. 13, 613--663. 
\bibitem{hori} K. Hori. \textit{Constraints For Topological Strings In $D\geq 1$. } Nucl. Phys. B439 (1995) 395--420.
\bibitem{Jin1} M. Jinzenji. \textit{Mirror Map as Generating Function of Intersection Numbers: Toric Manifolds with Two K\"{a}hler Forms}, arXiv:1006.0607, Comm. Math. Phys. 323(2013), no. 2, 747--811. 
\bibitem{Jin3}M. Jinzenji. \textit{Classical Mirror Symmetry.} SpringerBriefs in Mathematical Physics, 29. Springer, Singapore, 2018. viii+140 pp. ISBN: 978-981-13-0055-4; 978-981-13-0056-1 81-01.

\bibitem{JM} M. Jinzenji, K. Matsuzaka. \textit{Period Integrals (Govental's $I$-function) of Calabi-Yau Hypersurface in $CP^{N-1}$ and Intersection Numbers of Moduli Space of Quasimaps from $CP^{1}$with Two Marked Points to $CP^{N-1}$.} arXiv:2206.06591, Preprint.
\bibitem{JS} M. Jinzenji, M. Shimizu, \textit{Multi-Point Virtual Structure Constants and Mirror Computation of  $CP^{2}$-model.} Commun.Num.Theor Phys. 07 (2013) 411--468. 
\bibitem{LM} A. Losev, Y. Manin.  \textit{New moduli spaces of pointed curves and pencils of flat connections. Dedicated to William Fulton on the occasion of his 60th birthday.} 
Michigan Math. J. 48 (2000), 443--472.
\bibitem{MOP} A. Marian, D. Oprea and R. Pandharipande, \textit{The moduli space of stable quotients}, Geom. Topol. \textbf{15}(3) (2011), 1651--1706.
\bibitem{S} H. Saito. \textit{Chow Rings of $\widetilde{Mp}_{0,2}(N,d)$ and $\overline{M}_{0,2}({\mathbb P}^{N-1},d)$ and Gromov-Witten Invariants of Projective Hypersurfaces of Degree 1 and 2}, Internat. J. Math. 28(2017), no. 12, 1750090.
\bibitem{zinger}A. Zinger.\textit{ The reduced genus-one  Gromov-Witten invariants of Calabi-Yau hypersurfaces.} J. Amer. Math. Soc. 22 (2009), 691-737.
\end{thebibliography}
\end{document}